\documentclass[11pt]{article}
\usepackage{amsmath, amsthm, amssymb, amsfonts, mathrsfs, mathtools}
\usepackage{graphicx}
\usepackage{makeidx}
\usepackage[left=2cm,right=2cm,top=2cm,bottom=2cm]{geometry}
\usepackage{caption}
\usepackage{subcaption}
\usepackage[section]{placeins}
\usepackage{enumitem}
\usepackage{tikz}
\usepackage{stmaryrd}

 \usetikzlibrary{decorations.pathreplacing}
\tikzstyle{vertex}=[auto=left,circle,draw=black,fill=white, inner sep=1.5]
\usepackage{hyperref}
\hypersetup{colorlinks=true, citecolor={blue}, linkcolor=blue, filecolor=magenta, urlcolor=cyan}

\newtheorem{theorem}{Theorem}[section]
\newtheorem{prop}[theorem]{Proposition}
\newtheorem{lema}[theorem]{Lemma}
\newtheorem{corollary}[theorem]{Corollary}

\theoremstyle{definition}

\newcommand{\Cay}{\operatorname{Cay}}

\title{Quadratic unitary Cayley graphs over finite commutative rings}

\author{Priya Dhankhar and Sanjay Kumar Singh\\
Department of Mathematics\\
Indian Institute of Science Education and Research Bhopal, India.\\
priya22@iiserb.ac.in\\ sanjayks@iiserb.ac.in}

\date{}

\begin{document}
	\maketitle
	
	\vspace{-0.3in}

\begin{abstract}
Let $R$ be a finite commutative ring
with identity, and let
$T_R=Q_R\cup(-Q_R)$, where
$Q_R=\{u^2:u\in R^\times\}$. 
In this paper, we investigate the spectral properties of
quadratic unitary Cayley graphs $\operatorname{Cay}(R,T_R)$, over finite commutative rings with identity. We determine all the eigenvalues of $\operatorname{Cay}(R,T_R)$
when $R$ has odd characteristic.
Our results extend the previously known spectral formulas for this family of graphs. We also obtain all the
eigenvalues when the characteristic of $R$ is
$2$. Furthermore, for rings of odd characteristic, we give a complete characterization of when $\operatorname{Cay}(R,T_R)$ is a strongly regular Cayley graph.
\end{abstract}


\section{Introduction}
A \textit{graph} $G$ is an ordered pair $(V(G), E(G))$ consisting of a non-empty vertex set $V(G)$ and an edge set $E(G)$ of unordered pairs of elements of $V(G)$. A graph is finite if $V(G)$ and $E(G)$ are finite sets. In this paper, we consider only finite graphs. The \textit{$(0,1)$-adjacency matrix} of $G$, denoted by $\mathcal{A}(G)$, is the square matrix $[a_{uv}]$, where $a_{uv}$ is given by 
\[a_{uv} = \left\{ \begin{array}{rl}
	1 &\mbox{ if }
	(u,v)\in E(G) \\ 
	0 &\textnormal{ otherwise.}
\end{array}\right. \] The eigenvalues of $\mathcal{A}(G)$ are called the \textit{eigenvalues} of $G$. If $\lambda_1, \ldots ,\lambda_k$ are the eigenvalues of a graph $G$ with multiplicities $m_1,\ldots ,m_k$, respectively, then 
$\text{Spec}(G)=\begin{pmatrix}
\lambda_1& \ldots& \lambda_k\\
m_1 & \ldots & m_k
\end{pmatrix}$ describes the spectrum of $G$.

In general, a Cayley graph is defined on a group by taking the connection set as any subset of a group. However, in this paper, we consider Cayley graphs only on finite commutative rings with a particular type of connection set. Let \(R\) be a finite commutative ring with unity and let
\(R^{\times}\) be the set of units of the ring \(R\). A well-known family of Cayley graphs is the family of unitary Cayley graphs, defined by $\operatorname{Cay}(R,R^{\times})$. D. Kiani et al.  \cite{kiani2011energy} computed the eigenvalues of \emph{unitary Cayley graphs} $\operatorname{Cay}(R,R^{\times})$ for any finite commutative ring $R$ with unity. Later on, Priya and Sanjay Kumar Singh~\cite{kumar2024spectral} calculated the eigenvalue of the Cayley graph ${\rm Cay}(R, xR^{\times})$ for a finite commutative ring $R$ with unity and $x\in R$.
Define  
\(
Q_R=\{r^2 \mid r\in R^{\times}\}
\)
and
\[
T_R=Q_R\cup (-Q_R)
\]
(not necessarily disjoint). The Cayley graph of $R$ with the connection set $T_R$, denoted by ${\rm Cay}(R,T_R)$, is a graph with $V({\rm Cay}(R,T_R))=R$ and $$E({\rm Cay}(R,T_R))=\{ (u,v)\colon u,v\in R, u-v \in T_R \}.$$ The Cayley graph ${\rm Cay}(R, T_R)$ is known as a \emph{quadratic unitary Cayley graph}.

Quadratic unitary Cayley graphs were introduced by Beaudrap \cite{Beaudrap2010Quadratic} in 2010 in the setting of the ring
\(\mathbb Z_n\), where he studied the Cayley graph $\operatorname{Cay}(\mathbb Z_n,T_{\mathbb Z_n})$. Later, in 2015, Liu and Zhou~\cite{Liu.Zhou2015Quadratic}
extended this notion from $\mathbb{Z}_n$ to arbitrary finite
commutative rings and studied several spectral properties of
these graphs under certain additional assumptions.
More precisely, let $R$ be a finite commutative ring with the
decomposition
\[
R=R_1\times R_2\times\cdots\times R_s,
\]
where each $R_i$ is a local ring with maximal ideal $M_i$.
Liu and Zhou~\cite{Liu.Zhou2015Quadratic} determined the spectrum
of the quadratic unitary Cayley graph ${\rm Cay}(R,T_R)$
under the assumption that
$|R_i|/|M_i|\equiv 1\pmod{4}$ for all $1\leq i\leq s$.
They also determined the spectrum of
${\rm Cay}(R_0\times R,T_{R_0\times R})$,
where $R_0$ is a finite local ring with maximal ideal $M_0$
satisfying $|R_0|/|M_0|\equiv -1\pmod{4}$.  Later on, in 2022, Huang \cite{Jing_Huang2022Quadratic} determined the spectrum of the quadratic unitary Cayley graph \(\rm{Cay}(\mathbb{Z}_n,T_{\mathbb{Z}_n})\) for any $n\geq 1$.

However, the spectral results of Liu and Zhou do not
cover the general case in which two or more residue
fields have orders congruent to $3$ modulo $4$.
Moreover, Huang's results concern the particular
class of rings $\mathbb{Z}_n$.

Motivated by these results, we study quadratic unitary
Cayley graphs over arbitrary finite commutative rings.
For rings of odd characteristic, we determine the
eigenvalues without imposing any restriction on the
number of residue fields whose orders are congruent
to $3$ modulo $4$ (see Theorem~\ref{spec_of_odd charc_ring}). Our approach uses additive
characters and quadratic Gauss sums over finite fields.

We also determine the spectrum of quadratic unitary
Cayley graphs over finite commutative rings of
characteristic $2$ (Theorem~\ref{spec_of_2k_char_ring}). Finally, using the spectral
results obtained for rings of odd characteristic,
we give a complete characterization of the rings of odd characteristic
for which $\Cay(R,T_R)$ is strongly regular (Theorem~\ref{srg_charc_for_odd_ring}).

This paper is organized as follows.
Section~2 contains some basic definitions and preliminary results.
In Sections~3 and~4, we determine the spectrum of quadratic unitary
Cayley graphs over finite commutative rings of odd characteristic
and characteristic $2$,
respectively.
Finally, in Section~5, we give a complete characterization of
strongly regular quadratic unitary Cayley graphs over
finite commutative rings of odd characteristic.

\section{Preliminaries}
In this section, we recall some basic definitions, notation,
and preliminary results concerning finite commutative rings, and additive characters that will be used
throughout the paper.
Let \(R\) be a finite commutative ring with unity. The \emph{Jacobson radical} of \(R\), denoted by \(J(R)\), is defined as the intersection of all maximal ideals of \(R\). It is well known that
\[
R^\times/(1+J(R))\cong (R/J(R))^\times
\]
(see \cite[Lemma~2.2]{stewart1972finite}) for more details.
\begin{lema}\label{R*_in terms_of_JR}
     Let $R$ be a finite commutative ring and $M_1,\ldots,M_k$ be all distinct maximal ideals of \(R\). Then \[|R^\times|
=\left|J(R)\right|
\prod_{i=1}^{k}
\left(\frac{|R|}{|M_i|}-1\right).\]
\end{lema}
\begin{proof}
By the Chinese Remainder Theorem, we have $R/J(R)\cong \prod_{i=1}^{k}R/M_i$.
Since $R^\times/(1+J(R))\cong (R/J(R))^\times$, we obtain
\[
|R^\times|
=|J(R)|\prod_{i=1}^{k}|(R/M_i)^\times|
=|J(R)|\prod_{i=1}^{k}
\left(\frac{|R|}{|M_i|}-1\right).
\]
This completes the proof.
\end{proof}

\begin{lema}\label{square_map_bij_on_1+J(R)}
Let $R$ be a finite commutative ring with unity
and of odd characteristic. Then the squaring map $\phi:1+J(R)\longrightarrow 1+J(R)$, defined as $\phi(u)=u^2$ is a bijection.
\end{lema}
\begin{proof}
Since $R$ has odd characteristic, $1+J(R)$ is a
finite abelian group of odd order. Let
$|1+J(R)|=2m+1$ for some nonnegative integer $m$. By Lagrange's theorem, for every $u\in 1+J(R)$,
we have $u^{2m+1}=1$.
Consequently, $(u^{m+1})^2=u^{2m+2}=u$.
Thus, every element of $1+J(R)$ is a square
of an element of $1+J(R)$. Hence, the squaring map
is surjective. Since $1+J(R)$ is finite,
$\phi$ is also injective and hence bijective.
\end{proof}

\begin{prop}\label{epression_of_T_r when R/M=1MOD4 AND 3 MOD 4}
Let $R$ be a finite commutative ring with identity, and let
$M$ be a maximal ideal of $R$. Suppose that $q=\frac{|R|}{|M|}$
is odd.
Then
\[
T_{R/M}=
\begin{cases}
Q_{R/M}, & \text{if } q\equiv 1\pmod{4},\\[4pt]
(R/M)^\times, & \text{if } q\equiv 3\pmod{4}.
\end{cases}
\]
\end{prop}
\begin{proof}
Since $R/M$
is a finite field of order $q$. Therefore,
$(R/M)^\times$ is a cyclic group of order $q-1$.
Let $g$ be a generator of $(R/M)^\times$. Then $Q_{R/M}=\langle g^2\rangle$ and $|Q_{R/M}|=\frac{q-1}{2}$.
First, suppose that $q\equiv 1\pmod{4}$.
Then $4\mid(q-1)$, and hence 
\[
-1=g^{(q-1)/2}
=\left(g^{(q-1)/4}\right)^2
\in Q_{R/M}.
\]
Thus, $-Q_{R/M}=Q_{R/M}$, which implies $T_{R/M}=Q_{R/M}$. Next, suppose that $q\equiv 3\pmod{4}$.
Then $(q-1)/2$ is odd, so $-1=g^{(q-1)/2}$ does not belong to $Q_{R/M}$. Consequently, $Q_{R/M}\cap(-Q_{R/M})=\varnothing$.
Since $|Q_{R/M}|=|-Q_{R/M}|=\frac{q-1}{2}$,
we obtain $|T_{R/M}|
=|Q_{R/M}|+|-Q_{R/M}|
=q-1$.
As $T_{R/M}\subseteq(R/M)^\times$ and
$|(R/M)^\times|=q-1$, it follows that $T_{R/M}=(R/M)^\times$.
This completes the proof.
\end{proof}
An \emph{additive character} of $R$ is a group
homomorphism $\psi:(R,+)\longrightarrow \mathbb{C}^{\times}$.
Thus,
\[
\psi(x+y)=\psi(x)\psi(y)
\quad\text{for all }x,y\in R.
\]
To determine the eigenvalues of quadratic unitary Cayley graphs, we use additive characters and quadratic Gauss sums over finite fields. We first recall some basic definitions and results concerning these characters. 

An \emph{additive character} of $R$ is a group
homomorphism $\psi:(R,+)\longrightarrow \mathbb{C}^{\times}$.
Thus,
\[
\psi(x+y)=\psi(x)\psi(y)
\quad\text{for all }x,y\in R.
\]An additive character $\psi$ is called trivial if
$\psi(x)=1$ for every $x\in R$.Since the additive group of $R$ is finite
and abelian, its additive characters can be
indexed by the elements of $R$. We fix such
an indexing and denote these characters by
$\psi_\alpha$, where $\alpha\in R$. We choose this indexing
to be compatible with the parametrization of
the additive characters of $R/J(R)$ given below.

 In 1979, Babai~\cite{babai1979spectra}, determined the
eigenvalues of Cayley graphs using characters.
In particular, his result implies that the
eigenvalues of $\Cay(R,T_R)$ are given by
\[
\lambda_\alpha=\sum_{s\in T_R}\psi_\alpha(s).
\]
\begin{itemize}
\item Let $F$ be a finite field of odd order $q$.
The \emph{quadratic character} of $F$ is the
multiplicative character
$\eta:F^\times\longrightarrow\{1,-1\}$
defined by
\[
\eta(x)=
\begin{cases}
1, & \text{if }x\in Q_F,\\
-1, & \text{if }x\in F^\times\setminus Q_F,
\end{cases}
\]
We extend $\eta$ to $F$ by setting $\eta(0)=0$.
\item 
Let $F$ be a finite field of odd order $q$,
and let $\eta$ denote the quadratic character
of $F$. Fix a nontrivial additive character
$\theta:F\to\mathbb{C}^{\times}$ such that
\[
\sum_{x\in F^\times}\eta(x)\theta(x)
=
\begin{cases}
\sqrt{q}, & q\equiv 1\pmod{4},\\
i\sqrt{q}, & q\equiv 3\pmod{4}.
\end{cases}
\]
For each $a\in F$, define $\psi_a(x)=\theta(ax)$ for all $x\in F$.
Then every additive character of $F$ is of
the form $\psi_a$ for a unique $a\in F$.
Moreover, for $a\ne 0$, we have
\[
\sum_{x\in F^\times}\eta(x)\psi_a(x)
=
\eta(a)\sum_{x\in F^\times}\eta(x)\theta(x).
\]

In particular, when $F=R/M$, we use the
notation $\chi_{a+M}$ for the additive
character defined by this convention.
\item Every additive character $\psi_{\alpha}$ of $R$
that is trivial on $J(R)$ induces a unique
additive character $\overline{\psi_{\alpha}}$ of $R/J(R)$,
defined by
\[
\overline{\psi_{\alpha}}(r+J(R))=\psi_{\alpha}(r),
\qquad r\in R.
\]
Conversely, every additive character of
$R/J(R)$ induces an additive character of $R$
that is trivial on $J(R)$ by composition
with the natural projection.
\item By the Chinese Remainder Theorem, we have
\[
R/J(R)\cong\prod_{i=1}^{k}R/M_i.
\]
For each $i$, fix a normalized nontrivial
additive character $\theta_i$ of $R/M_i$.
For $\alpha\in R$, define
\[
\chi_{\alpha+M_i}(x+M_i)
=
\theta_i((\alpha+M_i)(x+M_i)).
\]
Then every additive character of $R/J(R)$
is uniquely of the form
\[
\overline{\psi}_{\alpha}(x+J(R))
=
\prod_{i=1}^{k}
\chi_{\alpha+M_i}(x+M_i),
\]
where $\alpha+J(R)\in R/J(R)$.
\item Let $S$ be a fixed set of representatives
of the cosets of $J(R)$ in $R$.
We choose the indexing of the additive characters
of $R$ such that, for every $\alpha\in S$,
\[
\psi_\alpha(x)
=
\overline{\psi}_\alpha(x+J(R))
=
\prod_{i=1}^{k}
\chi_{\alpha+M_i}(x+M_i),
\qquad x\in R.
\]
The remaining additive characters of $R$
are indexed by the elements of $R\setminus S$.
With this choice of indexing, $J(R)\subseteq\ker(\psi_\alpha)
\quad\Longleftrightarrow\quad
\alpha\in S$.
\end{itemize}

\begin{lema}\label{eigenvaluesofpayleygraph}
Let $R$ be a finite commutative ring with unity, and let
$M$ be a maximal ideal of $R$ such that $q=|R/M|$ is odd.
Let $\chi_{a+M}$ be an additive character of $R/M$. Put
\[
G_q=
\begin{cases}
\sqrt{q}, & \text{if } q\equiv 1\pmod{4},\\[4pt]
i\sqrt{q}, & \text{if } q\equiv 3\pmod{4}.
\end{cases}
\]
Then
\[
\sum_{x\in Q_{R/M}}\chi_{a+M}(x)
=
\begin{cases}
\dfrac{q-1}{2},
& \text{if } a+M=M,\\[8pt]
\dfrac{-1+G_q}{2},
& \text{if } a+M\in Q_{R/M},\\[8pt]
\dfrac{-1-G_q}{2},
& \text{if } a+M\in
(R/M)^\times\setminus Q_{R/M}.
\end{cases}
\]
\end{lema}

\begin{proof}
Let $\eta$ denote the quadratic character of $R/M$.
If $a+M=M$, then $\chi_{a+M}$ is the trivial
additive character. Therefore, $\sum_{x\in Q_{R/M}}\chi_{a+M}(x)
=|Q_{R/M}|=\frac{q-1}{2}$.
Now, suppose that $a+M\neq M$.
Since
\[
\frac{1+\eta(x)}{2}
=
\begin{cases}
1, & \text{if } x\in Q_{R/M},\\[4pt]
0, & \text{if } x\in (R/M)^\times\setminus Q_{R/M}.
\end{cases}
\]
we obtain
\[
\begin{aligned}
\sum_{x\in Q_{R/M}}\chi_{a+M}(x)
&=\frac12\sum_{x\in(R/M)^\times}
(1+\eta(x))\chi_{a+M}(x)\\
&=\frac12\left(
\sum_{x\in(R/M)^\times}\chi_{a+M}(x)
+
\sum_{x\in(R/M)^\times}
\eta(x)\chi_{a+M}(x)
\right)\\
&=\frac12\left(
-1+\eta(a+M)
\sum_{x\in(R/M)^\times}\eta(x)\theta(x)
\right)\\
&=\frac12\left(-1+\eta(a+M)G_q\right).
\end{aligned}
\]
Since $\eta(a+M)=1$ when $a+M\in Q_{R/M}$
and $\eta(a+M)=-1$ otherwise. This completes the proof.
\end{proof}

\section{Spectrum of quadratic unitary Cayley graphs over rings of odd characteristic}

In this section, we determine the complete spectrum of quadratic
unitary Cayley graphs over finite commutative rings of odd
characteristic. We first establish some preliminary results
concerning character sums, which will be used to obtain
the eigenvalues of $\operatorname{Cay}(R,T_R)$.


\begin{lema}\label{T_R as disjoint union of costs of JR}
Let $R$ be a finite commutative ring with identity and of odd characteristic. Then there exists a set $\mathcal{C}\subseteq T_R$ such that
\[
T_R=\bigcup_{a\in\mathcal{C}}\bigl(a+J(R)\bigr),
\]
where the cosets $a+J(R)$, $a\in\mathcal{C}$, are pairwise distinct.
\end{lema}

\begin{proof}
We first show that $Q_R$ is a union of cosets of $J(R)$. Let $a\in Q_R$. Then $a=u^2$ for some $u\in R^\times$. For any $x\in J(R)$, $a+x = u^2+x = u^2\left(1+u^{-2}x\right)$.
Since $J(R)$ is an ideal, $u^{-2}x\in J(R)$. Therefore, $1+u^{-2}x\in 1+J(R)$. By Lemma \ref{square_map_bij_on_1+J(R)}, there exists $v\in 1+J(R)$ such that $1+u^{-2}x=v^2$. Hence, $a+x=u^2v^2=(uv)^2\in Q_R$.
Thus, $a+J(R)\subseteq Q_R$ for every $a\in Q_R$. Since $a\in a+J(R)$, it follows that
\[
Q_R=\bigcup_{a\in Q_R}\bigl(a+J(R)\bigr).
\]
Hence, $Q_R$ is a union of cosets of $J(R)$.
Similarly, we can show that $-Q_R$ is also a union of cosets of $J(R)$. Consequently, $T_R=Q_R\cup(-Q_R)$ is a union of cosets of $J(R)$. After removing repeated cosets we obtain a set of representatives $\mathcal{C}\subseteq T_R$ such that
\[
T_R=\bigcup_{a\in\mathcal{C}}\bigl(a+J(R)\bigr).
\]
Therefore, $T_R$ is a disjoint union of cosets of $J(R)$ in $R$.
\end{proof}

\begin{lema}\label{lambda_alpha_in_terms_of_JR and chractersum}
Let $R$ be a finite commutative ring with identity
and of odd characteristic. Then the eigenvalue
$\lambda_\alpha$ of $\operatorname{Cay}(R,T_R)$
satisfies
\[
\lambda_\alpha=
\begin{cases}
|J(R)|\displaystyle\sum_{a\in C}\psi_\alpha(a),
& \text{if }J(R)\subseteq\langle\alpha\rangle^\perp,\\[6pt]
0, & \text{otherwise},
\end{cases}
\]
where $C\subseteq T_R$ is a set of representatives
of the distinct cosets of $J(R)$ contained in $T_R$.
\end{lema}

\begin{proof}
By Lemma \ref{T_R as disjoint union of costs of JR}, we have \[
T_R=\bigcup_{a\in C}(a+J(R)).
\]
Therefore,
\[
\begin{aligned}
\sum_{s\in T_R}\psi_\alpha(s)
&=\sum_{a\in C}\sum_{x\in J}
  \psi_\alpha(a+x)\\
&=\sum_{a\in C}\psi_\alpha(a)
  \sum_{x\in J}\psi_\alpha(x).
\end{aligned}
\]
By the orthogonality relation for additive
characters, we have
\[
\sum_{x\in J}\psi_\alpha(x)
=
\begin{cases}
|J|, & \text{if }
J(R)\subseteq\ker(\psi_\alpha),\\
0, & \text{otherwise}.
\end{cases}
\]
Hence, the result follows.
\end{proof}
In particular, $\lambda_{\alpha}=0$ whenever $J(R)\not\subseteq\ker(\psi_\alpha)$. Therefore, from now on, we only need to consider those characters \(\psi_{\alpha}\) which are trivial on \(J(R)\). Equivalently, we assume that $J(R)\subseteq\ker(\psi_\alpha)$.

\begin{lema}\label{expression_of_eigen_values_in _eigenvaluesof RJR}
Let $R$ be a finite commutative ring with identity
and of odd characteristic. Let $\psi_\alpha$ be an
additive character of $R$ such that
$J(R)\subseteq\ker(\psi_\alpha)$, and let
$\overline{\psi}_\alpha$ be the induced additive
character of $R/J(R)$. Then
\[
\sum_{s\in T_R}\psi_\alpha(s)
=
|J(R)|
\sum_{s\in T_{R/J(R)}}
\overline{\psi}_\alpha(s).
\]
\end{lema}

\begin{proof}
Let $\pi:R\longrightarrow R/J(R)$ be the
canonical projection. Since
$\pi(R^\times)=(R/J(R))^\times$, we have
\[
\pi(T_R)=T_{R/J(R)}.
\]
By Lemma~\ref{T_R as disjoint union of costs of JR}, $T_R$ is a disjoint union
of cosets of $J(R)$. Hence, each element
of $T_{R/J(R)}$ has exactly $|J(R)|$
preimages in $T_R$. Since $J(R)\subseteq\ker(\psi_\alpha)$,
the character $\psi_\alpha$ is constant
on each coset of $J(R)$. Therefore,
\[
\begin{aligned}
\sum_{s\in T_R}\psi_\alpha(s)
&=
\sum_{u\in T_{R/J(R)}}
\sum_{s\in\pi^{-1}(u)}
\psi_\alpha(s)\\
&=
|J(R)|
\sum_{u\in T_{R/J(R)}}
\overline{\psi}_\alpha(u).
\end{aligned}
\]
This completes the proof.
\end{proof}

Consequently, using the representation of additive
characters of $R/J(R)$ given in Section~2, we obtain
\[
\lambda_\alpha
=
|J(R)|
\sum_{(x_1,\ldots,x_k)\in T_{R/J(R)}}
\prod_{i=1}^{k}
\chi_{\alpha_i+M_i}(x_i),
\]
where $J(R)\subseteq\ker(\psi_\alpha)$.

\begin{lema}\label{TR_cardinality}
Let $R$ be a finite commutative ring of odd
characteristic, and let $M_1,\ldots,M_k$ be its
maximal ideals. Put $q_i=|R/M_i|$ for
$1\leq i\leq k$. Then, up to relabeling of the maximal ideals,
there exists an integer $s$ with $0\le s\leq k$
such that
\[
q_i \equiv
\begin{cases}
1 \pmod{4}, & 1\le i\le s,\\
3 \pmod{4}, & s+1\le i\le k.
\end{cases}
\]
If $s=k$, then $T_{R/J(R)}
=\prod_{i=1}^{k}T_{R/M_i}$.
If $0\leq s<k$, then
\[
T_{R/J(R)}
={}
\left(\prod_{i=1}^{s}T_{R/M_i}\right)
\times
\left(\prod_{i=s+1}^{k}Q_{R/M_i}\right)
\cup
\left(\prod_{i=1}^{s}T_{R/M_i}\right)
\times
\left(\prod_{i=s+1}^{k}(-Q_{R/M_i})\right).
\]
Moreover, the cardinality of $T_R$ is
\[
|T_R|=
\begin{cases}
\displaystyle
\frac{|J(R)|}{2^k}
\prod_{i=1}^{k}(q_i-1),
& s=k,\\[8pt]
\displaystyle
\frac{|J(R)|}{2^{k-1}}
\prod_{i=1}^{k}(q_i-1),
& 0\leq s<k.
\end{cases}
\]
\end{lema}

\begin{proof}
By the Chinese Remainder Theorem, $R/J(R)\cong\prod_{i=1}^{k}R/M_i$. Therefore, $(R/J(R))^{\times}\cong\prod_{i=1}^{k}(R/M_i)^\times$.
Since squaring is performed componentwise,
we have
\[
Q_{R/J(R)}
=\prod_{i=1}^{k}Q_{R/M_i} \hspace{2mm} \text{and} \hspace{2mm} -Q_{R/J(R)}
=\prod_{i=1}^{k}(-Q_{R/M_i}).
\]
By Proposition \ref{epression_of_T_r when R/M=1MOD4 AND 3 MOD 4}, if $q_i\equiv1\pmod4$,
then $-Q_{R/M_i}=Q_{R/M_i}=T_{R/M_i}$, whereas if $q_i\equiv3\pmod4$, then $Q_{R/M_i}\cap(-Q_{R/M_i})=\varnothing$.
Thus 
\[T_{R/J(R)}
={}
\left(\prod_{i=1}^{s}T_{R/M_i}\right)
\times
\left(\prod_{i=s+1}^{k}Q_{R/M_i}\right)
\cup
\left(\prod_{i=1}^{s}T_{R/M_i}\right)
\times
\left(\prod_{i=s+1}^{k}(-Q_{R/M_i})\right).\]
If $s=k$, then
\[
|T_{R/J(R)}|
=\prod_{i=1}^{k}|Q_{R/M_i}|
=\frac{1}{2^k}\prod_{i=1}^{k}(q_i-1).
\]
If $0\leq s<k$, the two sets in the
decomposition of $T_{R/J(R)}$ are disjoint.
Hence,
\[
\begin{aligned}
|T_{R/J(R)}|
&=2\prod_{i=1}^{k}|Q_{R/M_i}|
=\frac{1}{2^{k-1}}
\prod_{i=1}^{k}(q_i-1).
\end{aligned}
\]
By Lemma \ref{T_R as disjoint union of costs of JR}, $T_R$ is a union of cosets
of $J(R)$, and its image under the canonical
projection $R\to R/J(R)$ is $T_{R/J(R)}$.
Consequently, $|T_R|=|J(R)|\,|T_{R/J(R)}|$.
Hence, this completes the proof.
\end{proof}
\noindent

\begin{theorem} \label{spec_of_odd charc_ring}
Let $R$ be a finite commutative ring of odd
characteristic, and let $\{M_1,\ldots,M_k\}$
be the maximal ideals of $R$.
Put $q_i=|R/M_i|$ for each $1\le i\le k$.
Suppose that at least one $q_i$ is
congruent to $3$ modulo $4$.
Then, up to relabeling of the maximal ideals,
there exists an integer $s$ with $0\le s<k$
such that
\[
q_i \equiv
\begin{cases}
1 \pmod{4}, & 1\le i\le s,\\
3 \pmod{4}, & s+1\le i\le k.
\end{cases}
\]
Let \(\alpha\in R\), and let \(\psi_\alpha\) be an additive character such that $J(R)\subseteq\ker(\psi_\alpha).$ Put \[
A=\{i:\alpha+M_i\in Q_{R/M_i}\},
\]
and
\[
B=\{i:\alpha+M_i\in (R/M_i)^\times
\setminus Q_{R/M_i}\}.
\] Let $
D=(A\cup B)\cap\{s+1,\ldots,k\}$ and $|D|=m$. Then
\[
\begin{aligned}
\lambda_\alpha
={}&\frac{(-1)^{|B\cap\{1,\ldots,s\}|}|T_R|}
{\displaystyle\prod_{i\in D}(q_i-1)}
\prod_{i\in A\cap\{1,\ldots,s\}}
\frac{1}{\sqrt{q_i}+1}
\prod_{i\in B\cap\{1,\ldots,s\}}
\frac{1}{\sqrt{q_i}-1}\\
&\times\left[
\sum_{\substack{E\subseteq D\\
|E|\text{ even}}}
(-1)^{\,m-\frac{|E|}{2}+|E\cap B|}
\prod_{i\in E}\sqrt{q_i}
\right].
\end{aligned}
\]
\end{theorem}

\begin{proof}
By Lemma \ref{expression_of_eigen_values_in _eigenvaluesof RJR}, we have $\lambda_\alpha
=|J(R)|
\sum_{(x_1,\ldots,x_k)\in T_{R/J(R)}}
\prod_{i=1}^{k}\chi_{\alpha+M_i}(x_i).$
Using the decomposition of $T_{R/J(R)}$
given in Lemma \ref{TR_cardinality}, we obtain
\[
\begin{aligned}
\lambda_{\alpha}
=&
|J(R)|
\Bigg(
\prod_{i=1}^{s}
\sum_{x_i\in T_{R/M_i}}
\chi_{a+M_i}(x_i)
\prod_{i=s+1}^{k}
\sum_{x_i\in Q_{R/M_i}}
\chi_{a+M_i}(x_i)
\\
&+
\prod_{i=1}^{s}
\sum_{x_i\in \overline{T}_{R/M_i}}
\chi_{a+M_i}(x_i)
\prod_{i=s+1}^{k}
\sum_{x_i\in -Q_{R/M_i}}
\chi_{a+M_i}(x_i)
\Bigg)\\
=&2|J(R)|
\prod_{i=1}^{s}
\left(
\sum_{x_i\in T_{R/M_i}}
\chi_{\alpha+M_i}(x_i)
\right)
\operatorname{Re}\left[
\prod_{i=s+1}^{k}
\left(
\sum_{x_i\in Q_{R/M_i}}
\chi_{\alpha+M_i}(x_i)
\right)
\right].
\end{aligned}
\]

Using Lemma \ref{eigenvaluesofpayleygraph}, we have

\[
\begin{aligned}
\lambda_\alpha
={}&\frac{|J(R)|}{2^{k-1}}
\prod_{i\in(A\cup B)^c}(q_i-1)
\prod_{i\in A\cap\{1,\ldots,s\}}
(-1+\sqrt{q_i})
\prod_{i\in B\cap\{1,\ldots,s\}}
(-1-\sqrt{q_i})\\
&\times
\operatorname{Re}\Bigg[
\prod_{i\in A\cap\{s+1,\ldots,k\}}
(-1+\mathrm{i}\sqrt{q_i})
\prod_{i\in B\cap\{s+1,\ldots,k\}}
(-1-\mathrm{i}\sqrt{q_i})
\Bigg].
\end{aligned}
\]

By Lemma \ref{TR_cardinality}, $|T_R|=\frac{|J(R)|}{2^{k-1}}
\prod_{i=1}^{k}(q_i-1).$
Consequently,
\[
\begin{aligned}
\lambda_\alpha
={}&|T_R|
\prod_{i\in A\cap\{1,\ldots,s\}}
\frac{-1+\sqrt{q_i}}{q_i-1}
\prod_{i\in B\cap\{1,\ldots,s\}}
\frac{-1-\sqrt{q_i}}{q_i-1}\\
&\times
\operatorname{Re}\Bigg[
\prod_{i\in A\cap\{s+1,\ldots,k\}}
\frac{-1+\mathrm{i}\sqrt{q_i}}{q_i-1}
\prod_{i\in B\cap\{s+1,\ldots,k\}}
\frac{-1-\mathrm{i}\sqrt{q_i}}{q_i-1}
\Bigg].
\end{aligned}
\]
Since $\frac{-1+\sqrt{q_i}}{q_i-1}
=\frac{1}{\sqrt{q_i}+1}$
and $\frac{-1-\sqrt{q_i}}{q_i-1}
=-\frac{1}{\sqrt{q_i}-1}$, we have
\[
\begin{aligned}
\lambda_\alpha
={}&
\frac{(-1)^{|B\cap\{1,\ldots,s\}|}|T_R|}
{\displaystyle\prod_{i\in D}(q_i-1)}
\prod_{i\in A\cap\{1,\ldots,s\}}
\frac{1}{\sqrt{q_i}+1}
\prod_{i\in B\cap\{1,\ldots,s\}}
\frac{1}{\sqrt{q_i}-1}\\
&\times
\operatorname{Re}\left[
\prod_{i\in A\cap D}
(-1+\mathrm{i}\sqrt{q_i})
\prod_{i\in B\cap D}
(-1-\mathrm{i}\sqrt{q_i})
\right].
\end{aligned}
\]
It remains to evaluate the real part. For each
$i\in D$, put
\[
\varepsilon_i=
\begin{cases}
1, & i\in A\cap D,\\
-1, & i\in B\cap D.
\end{cases}
\]
Then
\[
\begin{aligned}
&\prod_{i\in A\cap\{s+1,\ldots,k\}}
(-1+\mathrm{i}\sqrt{q_i})
\prod_{i\in B\cap\{s+1,\ldots,k\}}
(-1-\mathrm{i}\sqrt{q_i})
=
\prod_{i\in D}
(-1+\varepsilon_i\mathrm{i}\sqrt{q_i}).
\end{aligned}
\]
Expanding this product over all subsets
$E\subseteq D$, we get
\[
\begin{aligned}
&\prod_{i\in D}
(-1+\varepsilon_i\mathrm{i}\sqrt{q_i})
=
\sum_{E\subseteq D}
(-1)^{m-|E|}
\mathrm{i}^{|E|}
\prod_{i\in E}\varepsilon_i\sqrt{q_i}.
\end{aligned}
\]
If $|E|$ is odd, then
$\mathrm{i}^{|E|}$ is purely imaginary.
Therefore, only subsets of even
cardinality contribute to the real part.
For even $|E|$, we have $\mathrm{i}^{|E|}
=(-1)^{|E|/2}.$
Moreover, $\prod_{i\in E}\varepsilon_i
=(-1)^{|E\cap B|}.$
Hence,
\[
\begin{aligned}
\operatorname{Re}\left[
\prod_{i\in A\cap\{s+1,\ldots,k\}}
(-1+\mathrm{i}\sqrt{q_i})
\prod_{i\in B\cap\{s+1,\ldots,k\}}
(-1-\mathrm{i}\sqrt{q_i})
\right]
=\sum_{\substack{E\subseteq D\\
|E|\text{ even}}}
(-1)^{m-|E|}
(-1)^{|E|/2}
(-1)^{|E\cap B|}
\prod_{i\in E}\sqrt{q_i}
\end{aligned}
\]
Substituting this into the preceding expression
for $\lambda_\alpha$, we obtain
\[
\begin{aligned}
\lambda_\alpha
={}&\frac{(-1)^{|B\cap\{1,\ldots,s\}|}|T_R|}
{\displaystyle\prod_{i\in D}(q_i-1)}
\prod_{i\in A\cap\{1,\ldots,s\}}
\frac{1}{\sqrt{q_i}+1}
\prod_{i\in B\cap\{1,\ldots,s\}}
\frac{1}{\sqrt{q_i}-1}\\
&\times\left[
\sum_{\substack{E\subseteq D\\
|E|\text{ even}}}
(-1)^{\,m-\frac{|E|}{2}+|E\cap B|}
\prod_{i\in E}\sqrt{q_i}
\right].
\end{aligned}
\]
This completes the proof.
\end{proof}

We now consider the case when $q_i\equiv 1\pmod{4}$
for all $1\leq i\leq k$. As a consequence of our results,
we recover the eigenvalue formula obtained by
Liu and Zhou~\cite{Liu.Zhou2015Quadratic}
in the following corollary.
\begin{corollary} \label{eigenvlues_odd_each_q_i_cong_1_mod_4}
    Let \(R\) be a finite commutative ring and \(M_1,M_2,\ldots,M_k\) be all the maximal ideals of \(R\) such that $q_i=|R/M_i|\equiv 1 \pmod{4}$ for all $1\leq i\leq k$. Let \(\psi_{\alpha}\) be an additive character of \(R\) such that $ J(R)\subseteq\ker(\psi_\alpha)$.
Then 
\[
\lambda_{\alpha}
=
(-1)^{|B|}
\frac{|R^\times|}
{
2^k
\displaystyle\prod_{i\in A}
\left(\sqrt{q_i}+1\right)
\displaystyle\prod_{j\in B}
\left(\sqrt{q_j}-1\right)
},
\]
 where $ A
=
\left\{
i\in\{1,\ldots,k\}:
\alpha+M_i\in Q_{R/M_i}
\right\}$ and $B
=
\left\{
i\in\{1,\ldots,k\}:
\alpha+M_i\in
\left(\frac{R}{M_i}\right)^{\times}
\setminus Q_{R/M_i}
\right\}$.
\end{corollary}
\begin{proof}
   We have $\lambda_{\alpha}= |J(R)|
\prod_{i=1}^{k}
\left(
\sum_{x_{i}\in T_{R/M_{i}}}
\chi_{\alpha+M_i}(x_{i})
\right)$. Using Lemma \ref{eigenvaluesofpayleygraph}, we have  

\begin{align*}
\lambda_\alpha &=  |J(R)|~~ \prod\limits_{i\not\in A\cup B} \frac{1}{2} \left( q_i-1 \right)   \prod\limits_{i\in A} \frac{1}{2} \left( \sqrt{q_i}-1 \right)  ~ \prod\limits_{i\in B} \frac{-1}{2}  \left( \sqrt{q_i} + 1 \right)       \\
&=
(-1)^{|B|}\frac{|J(R)|}
{2^k}
\frac{
\displaystyle
\prod_{i\notin A\cup B}(q_i-1)
\prod_{i\in A}(\sqrt{q_i}-1)(\sqrt{q_i}+1)
\prod_{i\in B}(\sqrt{q_i}+1)(\sqrt{q_i}-1)
}{
\displaystyle
\prod_{i\in A}(\sqrt{q_i}+1)
\prod_{i\in B}(\sqrt{q_i}-1)
}
\\[8pt]
&=
(-1)^{|B|}\frac{|J(R)|}
{2^k}
\frac{
\displaystyle\prod_{i=1}^{k}(q_i-1)
}{
\displaystyle
\prod_{i\in A}(\sqrt{q_i}+1)
\prod_{i\in B}(\sqrt{q_i}-1)
}.
\end{align*}
Using Lemma \ref{R*_in terms_of_JR}, we have
\begin{align*}
    \lambda_\alpha
&=
(-1)^{|B|}\frac{|R^\times|}
{
2^k
\displaystyle\prod_{i\in A}(\sqrt{q_i}+1)
\displaystyle\prod_{i\in B}(\sqrt{q_i}-1)
}
\end{align*}
This completes the proof.
\end{proof}


\section{Spectrum of quadratic unitary Cayley graphs over rings of characteristic \(2\)}

In this section, we determine the spectrum of quadratic unitary
Cayley graphs over finite commutative rings of characteristic $2$.
We begin by introducing some notation and recalling a few basic
observations that will be used throughout this section.

\begin{itemize}

\item Let $R$ be a finite commutative ring with unity and
of characteristic $2$. Then
\[
R \cong R_1 \times \cdots \times R_s,
\]
where each $R_i$ is a finite commutative local ring with
maximal ideal $M_i$ of order $m_i$, for each $1\leq i\leq s$.
By \cite[Theorem~8.7]{atiyah2018introduction},
this decomposition is unique up to the permutation of local rings.
\item Since $R=R_1\times\cdots\times R_s$, we have $Q_R=Q_{R_1}\times\cdots\times Q_{R_s}$. As the characteristic of $R$ is $2$, we have $-Q_R=Q_R$, and hence
$T_R=Q_R$. 
\item Let $R$ be a finite commutative local ring with
maximal ideal $M$. We define
\[
R^{[2]}=\{u^2:u\in R\}
\]
and
\[
M^{[2]}=\{m^2:m\in M\}.
\]

\item For a finite commutative ring
$R=R_1\times\cdots\times R_s$ of characteristic $2$,
we define
\[
R^{[2]}=\{x^2:x\in R\}.
\]
Clearly,
\[
R^{[2]}=R_1^{[2]}\times\cdots\times R_s^{[2]}.
\]

\end{itemize}
\begin{lema} \label{prop_of_local_ring_of_char_2}
Let $R$ be a finite commutative local ring of characteristic $2$.
Then $R^{[2]}$ is a subring of $R$, and $M^{[2]}$ is the unique
maximal ideal of $R^{[2]}$. Moreover,
\[
\frac{R^{[2]}}{M^{[2]}}
\cong
\frac{R}{M}, \qquad  \text{and}\qquad Q_R=\left(R^{[2]}\right)^\times.
\]
\end{lema}

\begin{proof}
Since $\operatorname{char}(R)=2$, we have $(a+b)^2=a^2+b^2$ and $(ab)^2=a^2b^2$
for all $a,b\in R$. Hence $R^{[2]}$ is a subring of $R$,
and $M^{[2]}$ is an ideal of $R^{[2]}$.
We first show that $R^{[2]}$ is a local ring with unique
maximal ideal $M^{[2]}$. Let $a^2\in R^{[2]}$.
If $a\notin M$, then $a\in R^\times$, since $R$ is local.
Consequently, $a^2(a^{-1})^2=1$, therefore $a^2$ is a unit of $R^{[2]}$.
Conversely, if $a^2$ is a unit of $R^{[2]}$, then
$a^2b^2=1$ for some $b\in R$. Thus $a(ab^2)=1$,
and hence $a\in R^\times$.
Therefore,
\[
(R^{[2]})^\times
=\{a^2:a\in R^\times\}
=Q_R.
\]
Since $R$ is local, an element $a\in R$ is a nonunit
if and only if $a\in M$. It follows that $R^{[2]}\setminus (R^{[2]})^\times=M^{[2]}$.
Hence $M^{[2]}$ is the unique maximal ideal of $R^{[2]}$.
Consider the ring homomorphism
\[
\phi:R^{[2]}\longrightarrow R/M,
\qquad
\phi(a^2)=a^2+M.
\]
Since $R/M$ is a finite field
of characteristic $2$, the Frobenius map
$x\mapsto x^2$ is bijective on $R/M$.
Hence $\phi$ is surjective.
Moreover,
\[
\ker(\phi)
=\{a^2\in R^{[2]}:a^2\in M\}.
\]
Since $M$ is a prime ideal, $a^2\in M$ if and only if
$a\in M$. Therefore, $\ker(\phi)=M^{[2]}$.
By the First Isomorphism Theorem, $R^{[2]}/M^{[2]}\cong R/M$.
This completes the proof.
\end{proof}

\begin{prop}\cite[Proposition 2.1]{kiani2011energy} \label{eigenvalues_unit_cayl_graph}
Let $R$ be a finite local ring with maximal ideal $M$ of size $m$. Then
\[
\operatorname{Spec}(\rm{Cay}(R,R^\times))
=
\left(
\begin{array}{ccc}
|R^\times| & -m & 0\\[1mm]
1 & \dfrac{|R|}{m}-1 & \dfrac{|R|}{m}(m-1)
\end{array}
\right)
=
\left(
\begin{array}{ccc}
|R^\times| & -m & 0\\[1mm]
1 & \dfrac{|R^\times|}{m} & \dfrac{|R|}{m}(m-1)
\end{array}
\right).
\]
\end{prop}

\begin{theorem}\label{eigenvalues_of_Quad_cay_graph_of_char_2_local}
Let $R$ be a finite commutative local ring of characteristic $2$ with
maximal ideal $M$ with $q=\left|\frac{R}{M}\right|$. Then 
\[
\operatorname{Spec}(\operatorname{Cay}(R,T_R))
=
\left\{
\begin{array}{ccc}
|Q_R|
&
-|M^{[2]}|
&
0
\\[2mm]
\displaystyle \frac{|R|}{|R^{[2]}|}
&
\displaystyle
\frac{|R|}{|R^{[2]}|}(q-1)
&
\displaystyle
\frac{|R|}{|R^{[2]}|}
\left(|R^{[2]}|-q\right)
\end{array}
\right\},
\]
\end{theorem}
\begin{proof}
By Lemma \ref{prop_of_local_ring_of_char_2}, $R^{[2]}$ is a finite local ring with unique maximal
ideal $M^{[2]}$, and
\[
Q_R=(R^{[2]})^\times.
\]
Since $\operatorname{char}(R)=2$, we have $-Q_R=Q_R$. Hence $T_R=Q_R$. As $Q_R\subseteq R^{[2]}$, two vertices of $R$ can be adjacent only if they belong to the same additive coset of $R^{[2]}$. Moreover, the subgraph induced by
each such coset is isomorphic to
\[
\operatorname{Cay}
\bigl(R^{[2]},(R^{[2]})^\times\bigr).
\] Thus $\operatorname{Cay}(R,Q_R)$ is the disjoint union of $\frac{|R|}{|R^{[2]}|}$
copies of the unitary Cayley graph $\operatorname{Cay}
\left(
R^{[2]},(R^{[2]})^\times
\right).$
From Proposition \ref{eigenvalues_unit_cayl_graph}, we have 
\[
\operatorname{Spec} (\operatorname{Cay}
\left(
R^{[2]},Q_R
\right))
=
\left(
\begin{array}{ccc}
|(R^{[2]})^\times|
&
-|M^{[2]}|
&
0
\\[2mm]
1
&
\dfrac{|R^{[2]}|}{|M^{[2]}|}-1
&
\dfrac{|R^{[2]}|}{|M^{[2]}|}
\left(|M^{[2]}|-1\right)
\end{array}
\right).
\]
By Lemma \ref{prop_of_local_ring_of_char_2}, $\frac{|R^{[2]}|}{|M^{[2]}|}
=
\left|\frac{R^{[2]}}{M^{[2]}}\right|
=
\left|\frac{R}{M}\right|
=q$.
Then \[
\operatorname{Spec} (\operatorname{Cay}
\left(
R^{[2]},Q_R
\right))
=
\left(
\begin{array}{ccc}
|(R^{[2]})^\times|
&
-|M^{[2]}|
&
0
\\[2mm]
1
&
q-1
&
\left(|R^{[2]}|-q\right)
\end{array}
\right).
\]
Since $\operatorname{Cay}(R,T_R)$ consists of
$|R|/|R^{[2]}|$ identical components, each multiplicity is multiplied by $|R|/|R^{[2]}|$. Hence,

\[
\operatorname{Spec}(\operatorname{Cay}(R,T_R))
=
\left(
\begin{array}{ccc}
|Q_R|
&
-|M^{[2]}|
&
0
\\[2mm]
\dfrac{|R|}{|R^{[2]}|}
&
\dfrac{|R|}{|R^{[2]}|}(q-1)
&
\dfrac{|R|}{|R^{[2]}|}
\left(|R^{[2]}|-q\right)
\end{array}
\right).
\]
\end{proof}
Let $G =(V(G),E(G))$ and $H=(V(H),E(H))$ be two graphs. Then the tensor product $G\otimes H $ is the graph with vertex set $V(G)\times V(H)$, and $((u,v),(u',v')) \in E(G\otimes H)$ if and only if $(u,u') \in E(G)$ and $(v,v') \in E(H)$. 

Now, we have the following theorem.
\begin{theorem}\label{iso tensor product of graph}
    Let $R=R_1\times R_2\times\cdots\times R_k$ be a finite commutative ring of characteristic $2$. Then \[
\operatorname{Cay}(R,T_R)
\cong
\bigotimes_{i=1}^k
\operatorname{Cay}(R_i,T_{R_i}).
\]
\end{theorem}
\begin{proof}
     As the characteristic of $R$ is $2$, $T_R=Q_R$. By our assumptions, we have $ R= R_1\times \cdots \times R_k$ and $Q_R=Q_{R_1}\times\cdots\times Q_{R_k}$. Note that, for any
$a=(a_1,\ldots,a_k)$ and $b=(b_1,\ldots,b_k)$
in $R$, $a-b\in Q_R$ if and only if $a_i-b_i\in Q_{R_i}$ for each $1\leq i\leq k$.
Hence, this completes the proof.
\end{proof}
 
\begin{theorem}\cite{west2001introduction}\label{Eigenvalues_of_tensor_product_of_graphs}
Let $G$ and $H$ be two graphs. If $\lambda_1, \ldots, \lambda_n$ are the eigenvalues of $G$ and $\mu_1,\ldots,\mu_m $ are the eigenvalues of $H$, then $\lambda_{i}\mu_{j}$,where $1\leq i\leq n$ and $1\leq j \leq m$ are the eigenvalues of $G\otimes H$. 
\end{theorem}

\begin{theorem} \label{spec_of_2k_char_ring}
Let $R=R_1\times R_2\times\cdots\times R_k$ be a finite commutative ring of characteristic $2$, where each $R_i$
is a local ring with maximal ideal $M_i$. Put $q_i=\left|\frac{R_i}{M_i}\right|$ for all $i=1,\ldots,k.$
Then the eigenvalues of the quadratic unitary Cayley graph
$\operatorname{Cay}(R,T_R)$ are 

\begin{enumerate}
\item[(i)] $\lambda_C
=
(-1)^{|C|}
\frac{|Q_R|}
{\displaystyle\prod_{j\in C}(q_j-1)}$
with multiplicity $\frac{|R|}{|R^{[2]}|}
\prod_{j\in C}(q_j-1)$ for every subset $C\subseteq\{1,2,\ldots,k\}$.

\item[(ii)] $0$ with multiplicity $|R|
-
\frac{|R|}
{\displaystyle\prod_{i=1}^k |M_i^{[2]}|}$.

\end{enumerate}
\end{theorem}
\begin{proof}
\begin{enumerate}[label=(\roman*)]
\item Let $C\subseteq\{1,\ldots,k\}$. By Theorem \ref{eigenvalues_of_Quad_cay_graph_of_char_2_local} and \ref{Eigenvalues_of_tensor_product_of_graphs}, $\lambda_C
=
(-1)^{|C|}
\prod_{j\in C}|M_j^{[2]}|
\prod_{j\notin C}|Q_{R_j}|$ is a nonzero eigenvalue with multiplicity $\prod_{i=1}^k\frac{|R_i|}{|R_i^{[2]}|}
\prod_{j\in C}(q_j-1)$.

Since $|Q_R|=\prod_{j=1}^s|Q_{R_j}|$ and  $R^{[2]}= \prod_{i=1}^s R_{i}^{[2]}$.
We have 
\[
\lambda_C
=
(-1)^{|C|}
\frac{|Q_R|}
{\displaystyle\prod_{j\in C}
\frac{|Q_{R_j}|}{|M_j^{[2]}|}}=\lambda_C
=
(-1)^{|C|}
\frac{|Q_R|}
{\displaystyle\prod_{j\in C}(q_j-1)}.
\]
with multiplicity
\[
\prod_{i=1}^k\frac{|R_i|}{|R_i^{[2]}|}
\prod_{j\in C}(q_j-1)
=
\frac{|R|}{|R^{[2]}|}
\prod_{j\in C}(q_j-1).
\]

\item Using
\[
\sum_{C\subseteq \{1,\ldots,k\}}
\prod_{j\in C}(q_j-1)
=
\prod_{i=1}^k\left(1+(q_i-1)\right)
=
\prod_{i=1}^k q_i.
\] We obtain 
\[ \frac{|R|}{|R^{[2]}|}\sum_{C\subseteq P}\left(\prod_{j\in C}(q_j-1)\right)
=
\frac{|R|}{|R^{[2]}|}
\prod_{i=1}^k q_i
=
\frac{|R|}
{\displaystyle\prod_{i=1}^k|M_i^{[2]}|}.
\]  is the multiplicity of all non zero eigenvalues.
Hence the multiplicity of the zero eigenvalue is
\[
|R|-
\frac{|R|}
{\displaystyle\prod_{i=1}^k|M_i^{[2]}|}.
\]
\end{enumerate}
\end{proof}

\section{Strongly regular quadratic unitary Cayley graphs over rings of odd characteristic}

In this section, we characterize the finite commutative rings
of odd characteristic for which the quadratic unitary Cayley
graph $\operatorname{Cay}(R,T_R)$ is strongly regular.
Using the spectral results obtained in Section~3,
we establish necessary and sufficient conditions for
$\operatorname{Cay}(R,T_R)$ to be strongly regular.

A Cayley graph is said to be a
$(v,r,\lambda,\mu)$-strongly regular graph
if the underlying group has order $v$,
the graph is $r$-regular, every pair of
adjacent vertices has exactly $\lambda$
common neighbors, and every pair of
nonadjacent vertices has exactly $\mu$
common neighbors. Throughout this section, we allow strongly regular graphs to be disconnected.

It is well known that a connected,
non-complete regular graph is strongly
regular if and only if its adjacency
matrix has exactly three distinct
eigenvalues. Moreover, a disconnected strongly
regular graph is a disjoint union of complete
graphs. The eigenvalues other than
$r$ are called the \emph{restricted
eigenvalues}.

In this section, we use the eigenvalues
obtained in Section~3 to characterize
the finite commutative rings of odd
characteristic for which the connection set $T_R$ gives a strongly regular cayley graph.

\begin{lema}\cite[Chapter 9]{brouwer2011spectra}\label{lem:complete-multipartite-srg}
Let $q,h\geq 2$ be integers. Then the complete
$r$-partite graph $K_{h,\ldots,h}$, with each
part of size $h$, is strongly regular with
parameters
\[
(v,d,\lambda,\mu)
=
\bigl(qh,(q-1)h,(q-2)h,(q-1)h\bigr).
\]
\end{lema}
\begin{theorem} \label{srg_charc_for_odd_ring}
Let $R$ be a finite commutative ring of odd characteristic,
and let $M_1,\ldots,M_k$ be its maximal ideals.
Put $q_i=|R/M_i|$ and $h=|J(R)|$.
Then $\operatorname{Cay}(R,T_R)$ is a 
strongly regular graph if and only if one of the
following holds:
\begin{enumerate}
\item $k=1$, $q_1\equiv3\pmod4$, and $h>1$;
\item $R\cong\mathbb F_q$, where $q\equiv1\pmod4$;
\item $R\cong\mathbb F_q\times\mathbb F_q$,
where $q\equiv3\pmod4$ and $q>3$;
\item $R\cong\mathbb F_3^k$, where $k\ge2$.
\end{enumerate}
\end{theorem}
\begin{proof}
Put $d=|T_R|$ and $h=|J(R)|$. 
We divide the proof into the following cases.

\medskip
\noindent
\textbf{Case 1: $q_i\equiv 1\pmod 4$ for some $i$.}

By Theorem~\ref{spec_of_odd charc_ring}, the additive characters that are
nontrivial only in the $i$-th coordinate give the
two eigenvalues
\[
\lambda_1=\frac{d}{\sqrt{q_i}+1},
\qquad
\lambda_2=-\frac{d}{\sqrt{q_i}-1}.
\]
Since $q_i\geq 5$, we have $\lambda_2<0<\lambda_1<d$.
Thus $\lambda_1$ and $\lambda_2$ are two distinct restricted eigenvalues. We first show that $k=1$. Suppose, to the contrary, that $k\geq 2$, and
choose $j\neq i$. We consider the following two possibilities.

\smallskip
\noindent
\emph{Subcase 1.1: $q_j\equiv 1\pmod 4$.}
By Theorem~\ref{spec_of_odd charc_ring}, a character nontrivial in exactly
the $i$-th and $j$-th coordinates gives the
eigenvalue
\[
\lambda_3
=
\frac{d}
{(\sqrt{q_i}+1)(\sqrt{q_j}+1)}.
\]
Since $\sqrt{q_j}+1>1$, we obtain $0<\lambda_3<\lambda_1$. Consequently, $\lambda_2<0<\lambda_3<\lambda_1<d$. Therefore, the adjacency matrix has at least
three distinct eigenvalues other than $d$.
Hence the graph cannot be strongly regular.

\smallskip
\noindent
\emph{Subcase 1.2: $q_j\equiv 3\pmod 4$.} By Theorem~\ref{spec_of_odd charc_ring}, a character nontrivial in the
$i$-th and $j$-th coordinates gives the eigenvalue
\[
\lambda_3
=
-\frac{d}
{(\sqrt{q_i}+1)(q_j-1)}.
\]
As $q_j\geq 3$, we have $-\frac{d}{\sqrt{q_i}-1}
<
-\frac{d}
{(\sqrt{q_i}+1)(q_j-1)}
<0.$
Thus $\lambda_2<\lambda_3<0<\lambda_1<d$.
Again, the adjacency matrix has at least three
distinct restricted eigenvalues.
Hence, in either case, $k\geq 2$ is impossible.
We conclude that $k=1$.
Next, suppose that $h>1$. 
By Lemma \ref{lambda_alpha_in_terms_of_JR and chractersum}, every additive character
$\psi_\alpha$ that is nontrivial on $J(R)$
gives the eigenvalue $\lambda_\alpha=0$.
In particular, if $J(R)\neq\{0\}$, then
such a character exists, and hence $0$
is an eigenvalue of $\operatorname{Cay}(R,T_R)$.
Since $\lambda_2<0<\lambda_1$, the eigenvalues  $\lambda_1,\quad \lambda_2,\quad 0$ are three distinct eigenvalues different from $d$. Therefore, the graph cannot be strongly regular. Consequently, $h=1$. Since $k=1$ and $J(R)=0$, 
we have
\[
R\cong\mathbb F_q,
\qquad q\equiv 1\pmod 4.
\]
Conversely, suppose that $R\cong\mathbb F_q$,
where $q\equiv 1\pmod 4$.
Then $\operatorname{Cay}(R,T_R)$ is the Paley graph
of order $q$. Its degree is $d=\frac{q-1}{2}$, and its two restricted eigenvalues are
\[
\frac{-1+\sqrt q}{2},
\qquad
\frac{-1-\sqrt q}{2}.
\]
These are distinct and different from $d$.
Therefore, the graph is strongly regular.\\
This proves assertion $(2)$.

\medskip
\noindent
\textbf{Case 2: $q_i\equiv 3\pmod 4$ for every $i$.}

We will consider the cases $k=1$ and $k\geq 2$.

\medskip
\noindent
\textbf{Case 2.1: $k=1$.} \\

Since $q_1\equiv 3\pmod 4$, by Proposition \ref{epression_of_T_r when R/M=1MOD4 AND 3 MOD 4}, we have $T_R=R^\times$. Two distinct vertices $x,y\in R$ are adjacent
if and only if $x-y\in R^\times$. Since $R$ is local, this is equivalent to $x-y\notin J(R)$. Therefore, $x$ and $y$ are adjacent precisely
when they belong to different cosets of $J(R)$. There are $q_1$ cosets of $J(R)$, each containing
$h$ elements. Hence $\operatorname{Cay}(R,T_R)
\cong K_{\underbrace{h,\ldots,h}_{q_1\text{ parts}}}$,
the complete $q_1$-partite graph with each part
of size $h$. If $h>1$, this graph is non-complete, and  by Lemma \ref{lem:complete-multipartite-srg}, this is 
strongly regular graph with parameters
\[
(v,d,\lambda,\mu)
=
\bigl(q_1h,(q_1-1)h,(q_1-2)h,(q_1-1)h\bigr).
\]
If $h=1$, the graph is the
complete graph $K_{q_1}$, which is excluded by our definition of strongly regular graph.
Therefore, when $k=1$, the graph is a strongly
regular if and only if $h>1$. This proves assertion $(1)$.

\medskip
\noindent
\textbf{Case 2.2: $k\geq 2$.}\\
By Theorem~\ref{spec_of_odd charc_ring}, the characters nontrivial in
exactly one coordinate give the eigenvalues
\[
\lambda_i=-\frac{d}{q_i-1},
\qquad 1\leq i\leq k.
\]
Moreover, for any two distinct indices
$i,j$, the characters nontrivial in exactly
the $i$-th and $j$-th coordinates give
\[
\lambda_{ij}^{\pm}
=
\frac{d(1\pm\sqrt{q_iq_j})}
{(q_i-1)(q_j-1)}.
\]
Notice that $\lambda_i<0,
\lambda_{ij}^{-}<0, \hspace{2mm} \text{and} \hspace{2mm}
\lambda_{ij}^{+}>0$.
We first show that all the residue-field
orders $q_i$ must be equal.
Suppose that $q_i\neq q_j$ for some $i\neq j$.
Then $-\frac{d}{q_i-1}
\neq
-\frac{d}{q_j-1}$.
Hence $\lambda_i$ and $\lambda_j$ are two
distinct negative eigenvalues. Furthermore, since $q_i,q_j\geq 3$ and
$q_i\neq q_j$, we have $0<\lambda_{ij}^{+}<d$.
Hence, the graph has at least three
distinct eigenvalues $\lambda_i,
\lambda_j,\hspace{2mm} \text{and } \hspace{2mm}
\lambda_{ij}^{+}$,
all different from $d$.
Therefore, 
\[
q_1=q_2=\cdots=q_k=q,
\qquad q\equiv 3\pmod 4.
\]
Substituting $q_i=q_j=q$ into the preceding
expressions, we obtain 
\[
\lambda_i=-\frac{d}{q-1},\qquad 
\lambda_{ij}^{+}
=\frac{d(q+1)}{(q-1)^2},
\qquad
\lambda_{ij}^{-}
=-\frac{d}{q-1}.
\]

Put
\[
r=\frac{d(q+1)}{(q-1)^2},
\qquad
s=-\frac{d}{q-1}.
\]
We now distinguish the cases $q>3$ and $q=3$.

\smallskip
\noindent
\emph{Subcase 2.2.1: $q>3$.}
Since $q\equiv 3\pmod 4$, we have $q\geq 7$.
It follows that $\lambda_i<0<\lambda_{ij}<d$.
We claim that $k=2$.
Suppose, to the contrary, that $k\geq 3$.
By Theorem~\ref{spec_of_odd charc_ring}, a character nontrivial in
exactly three coordinates, with the same choice
of quadratic sign in each coordinate, gives
the eigenvalue
\[
t=\frac{d(3q-1)}{(q-1)^3}.
\]
Since $q>3$, we have $t>0$. Moreover,
\[
\begin{aligned}
r-t
&=
\frac{d(q+1)}{(q-1)^2}
-
\frac{d(3q-1)}{(q-1)^3}\\
&=
\frac{dq(q-3)}{(q-1)^3}>0.
\end{aligned}
\]
Consequently, $s<0<t<r<d.$
Thus the adjacency matrix has at least three
distinct eigenvalues other than $d$,
contradicting strong regularity. Hence $k=2$. Next, we show that $h=1$.
By Lemma \ref{lambda_alpha_in_terms_of_JR and chractersum}, every additive character
$\psi_\alpha$ that is nontrivial on $J(R)$
gives the eigenvalue $\lambda_\alpha=0$. But $s<0<r<d$.
Therefore, $s,0,r$ are three distinct
eigenvalues other than $d$, which is impossible.
Hence $h=1$, so $J(R)=0$. Since $k=2$, it follows
that
\[
R\cong\mathbb F_q\times\mathbb F_q,
\qquad
q\equiv 3\pmod 4,\quad q>3.
\]
Conversely, let $R\cong\mathbb F_q\times\mathbb F_q,
\hspace{2mm} \text{with} \hspace{2mm} q>3,\hspace{2mm} q\equiv 3\pmod 4.$
By Theorem~\ref{spec_of_odd charc_ring}, the adjacency spectrum consists
of the degree $d=\frac{(q-1)^2}{2}$
and the two eigenvalues
\[
r=\frac{d(q+1)}{(q-1)^2}
=\frac{q+1}{2},
\qquad
s=-\frac{d}{q-1}
=-\frac{q-1}{2}.
\]
Since $s<0<r<d$, 
the degree $d$ has multiplicity one, and hence
the graph is connected. It is also non-complete and has exactly three
distinct adjacency eigenvalues. Therefore,
it is strongly regular.

This proves assertion $(3)$.

\smallskip
\noindent
\emph{Subcase 2.2.2: $q=3$.}\\
We now have $q_1=\cdots=q_k=3$ and $k\geq 2$.
First, suppose that $h=1$.
Since $J(R)=0$, we obtain $R\cong\mathbb F_3^k$.
The only nonzero square in $\mathbb F_3$ is $1$. Therefore, $T_R=
\bigl\{(1,\ldots,1),(-1,\ldots,-1)\bigr\}$. Let $u=(1,\ldots,1)\in\mathbb F_3^k$
and put $H=\langle u\rangle
=\{0,u,-u\}$.
Clearly, $H$ is a subgroup of order $3$.
For any $x\in\mathbb F_3^k$, its neighbors are
precisely $x+u$ and $x-u$. Moreover, $x$, $x+u$ and $x-u$ 
are pairwise adjacent vertices and form a triangle. There are no edges between distinct cosets
of $H$, since every element of $T_R$ belongs
to $H$.\\
Consequently,
\[
\operatorname{Cay}(R,T_R)
\cong
\underbrace{K_3\sqcup\cdots\sqcup K_3}
_{3^{k-1}\text{ copies}}.
\]
Since $k\geq 2$, this graph is noncomplete.
Every adjacent pair has exactly one common
neighbor, while every nonadjacent pair has
no common neighbor. Hence it is strongly regular with parameters
\[
(v,d,\lambda,\mu)
=
(3^k,2,1,0).
\]
Conversely, suppose that $h>1$.
We show that the graph cannot be strongly regular.
Consider quotient map
\[
R\longrightarrow R/J(R)\cong\mathbb F_3^k
\]
By Lemma \ref{T_R as disjoint union of costs of JR}, $T_R$ is a union of cosets of $J(R)$,
and its image under the quotient map
$\pi:R\to R/J(R)$ is $T_{R/J(R)}$.
Consequently, for any $x,y\in R$, we have
\[
\begin{aligned}
x\sim y
&\iff x-y\in T_R\\
&\iff (x-y)+J(R)\in T_{R/J(R)}\\
&\iff (x+J(R))\sim(y+J(R)).
\end{aligned}
\]
Thus, adjacency in $\operatorname{Cay}(R,T_R)$
depends only on the corresponding cosets of $J(R)$.
As shown above, the quotient graph is a
disjoint union of triangles. The preimage
of each vertex of a triangle is a coset of
$J(R)$ containing $h$ elements. Thus each triangle lifts to a complete tripartite graph $K_{h,h,h}$,
and consequently
\[
\operatorname{Cay}(R,T_R)
\cong
\underbrace{K_{h,h,h}\sqcup\cdots\sqcup K_{h,h,h}}
_{3^{k-1}\text{ copies}}.
\]
Choose two distinct vertices in the same part
of one component. They are nonadjacent and have
exactly $2h$ common neighbors. On the other hand, two vertices belonging to
different components are also nonadjacent, but
have no common neighbors.\\
Since $h>1$, both types of nonadjacent pairs
exist, and their numbers of common neighbors
are different. Therefore, the graph is not
strongly regular. Hence $h=1$ is necessary.\\
We conclude that, when $q=3$ and $k\geq 2$,
the graph is strongly regular if and only if
\[
R\cong\mathbb F_3^k.
\]
This proves assertion $(4)$.
\end{proof}

\section*{Acknowledgements} 

\noindent The first author is supported by a Senior Research Fellowship from CSIR, Government of India (File No. 09/1020(15619)/2022-EMR-I).

\noindent\textbf{Data Availability :}  No data was gathered or used in this paper, so a ``data availability statement'' is not applicable.

\noindent\textbf{Conflict of interest :} The author states that there is no conflicts of interest.



\end{document}